\documentclass[11pt,leqno]{article}
\usepackage{amsmath,amsthm}
\usepackage {latexsym}
\usepackage{amssymb}
\newcommand{\bea}{\begin{eqnarray}}\newcommand{\eea}{\end{eqnarray}}
\newcommand{\beq}{\begin{equation}}\newcommand{\ee}{\end{equation}}
\def\pa{\partial}

\newtheorem{theorem}{Theorem}[section]\newtheorem{lemma}[theorem]{Lemma}
\newtheorem{defi}[theorem]{Definition}
\newtheorem{prop}[theorem]{Proposition}

\theoremstyle{remark}\newtheorem{remark}[theorem]{Remark}

\def\bm{\left( \begin{array}{cc}}\def\endm{\end{array}\right)}
\newcommand{\eq}{\end{equation}}
\def\pa{\partial}
\def \rectangle#1#2{\hbox{\vrule\vbox to #2
              {\hrule\hbox to #1{\hfil}\vfil\hrule}\vrule}}

\def\ro{{\varrho_{\,}}}
\def\ringabove{\,{}^{{}^{{}^\circ}}\!\!\!\!}

\numberwithin{equation}{section}

\begin{document}
\author {Hans Lindblad\thanks{Part of this work was done while H.L. was a Member of the
Institute for  Advanced Study, Princeton, supported by the NSF
grant DMS-0111298 to the Institute. H.L. was also partially
supported by the NSF Grant DMS-0200226. } \, and Avy
Soffer\thanks{Also a member of the Institute of Advanced
Study, Princeton.Supported in part by NSF grant DMS-0100490.}\\
University of California at San Diego and Rutgers University}
\title {A Remark on Long Range Scattering for the  nonlinear Klein-Gordon
equation} \maketitle

\section{Introduction} We consider the problem of scattering
for the critical nonlinear Klein-Gordon in one space dimension:
 \beq \label{eq:nlkg}
 \square v+ v=-\beta v^3,
 \eq
 where $\square=\pa_t^2-\pa_x^2$.
 Recall first that a solution of
the linear Klein-Gordon, i.e. $\beta=0$, is asymptotically given
by (as $\ro$  tends to infinity)

 \beq \label{eq:linearasymptotic}
 u(t,x)\sim \ro^{-1/2} e^{i\ro}\, a(x/\!\ro)+
 \ro^{-1/2} e^{-i\ro}\, \overline{a(x/\!\ro)},
 \qquad \text{where}\quad \ro=(t^2-|x|^2)^{1/2}\geq 0.
 \eq
Here $a(x/\!\ro)=(t/\!\ro) \widehat{u}_+(-x/\!\ro)=\sqrt{1+x^2/\ro^2}\,
\widehat{u}_+(-x/\!\ro)$, where $\widehat{u}_+(\xi)=\int u_+(x) \,
e^{-ix \xi}\, dx$ denotes the Fourier transform with respect to
$x$ only, $\widehat{u}_+=(\widehat{u}_0-i(|\xi|^2+1)^{-1/2}
\widehat{u}_1)/2$, where $u_0=u\big|_{t=0}$ and $u_1=\pa_t
u\big|_{t=0}$. Here the right hand side is to be interpreted as
$0$ outside the light cone, when $|x|>t$.
\eqref{eq:linearasymptotic} can be proven using stationary phase,
see e.g. \cite{H1}, where a complete asymptotic expansion into
negative powers of $\ro$ was given. Recently, Delort\cite{D1}
proved that \eqref{eq:nlkg} with small initial data have a global
solution with asymptotics of the form
 \beq
\label{eq:nonlinearasymptotic}
 v(t,x)\sim \ro^{-1/2} e^{i\phi_0(\ro\!,\,x/\ro)}\, a(x/\!\ro)+
 \ro^{-1/2} e^{-i\phi_0(\ro\!,\, x/\ro)}\, \overline{a(x/\!\ro)},
 \qquad
 \phi_0(\ro,x/\!\ro)=\ro+\frac{3}{8}\beta |a(x/\!\ro)|^2\,\ln{\ro}
 \eq

 We consider the inverse problem of scattering, i.e. we show that
 for any given asymptotic expansion of the above form (1.3) there is a
 solution agreeing with it at infinity. More precisely, we show:
\begin{theorem}\label{theorem:exist} Suppose that $a$ and $b_1=b-b_0$ are fast decaying smooth real valued
functions, where $b_0$ is a constant and $|a^{(k)}(x/\!\ro)|
+|b_1^{(k)}(x/\!\ro)|\leq
C_{Nk j}(1+|x/\!\ro|)^{-N}$, for any $k\geq 0$ and $N$. Let \beq
\label{eq:nonlinearasymptotictwo}
 v_0(t,x)=\ro^{-1/2} \, a(x/\!\ro)\cos{\phi(\ro,x/\!\ro)}
 \qquad
 \phi(\ro,x/\!\ro)=\ro+\frac{3}{8} \beta\,
 |a(x/\!\ro)|^2\,\ln{\ro}+b(x/\!\ro),
 \eq
 interpreted as $0$ when $|x|\geq t$.
 Then  there is $T<\infty$ such that \eqref{eq:nlkg} has a smooth solution $v$ for $T\leq t<\infty$ satisfying $v\sim v_0$
as $t\to\infty$. More precisely for $t\geq T$ we have
\beq\label{eq:v0asym}
\| (v-v_0)(t,\cdot)\|_{L^\infty}+\sum_{|\alpha|\leq 1}\|\pa_{t,x}^\alpha(v-v_0)(t,\cdot)\|_{L^2}\leq
C(1+\ln (1+t))^2 (1+t)^{-1}.
\eq
\end{theorem}
\begin{remark} Note that the initial data $(a,b)$, as well as the parameter $\beta$, can be arbitrarily large.
\end{remark}
\begin{remark} We also get a
complete asymptotic expansion, see Theorem \ref{theorem:asymptotic}.
\end{remark}
\begin{remark} Note that by \eqref{eq:nonlinearasymptotictwo}-\eqref{eq:v0asym} the solution constructed has bounded energy
$$
\frac{1}{2}\int v_t(t,x)^2+v_x(t,x)^2+v(t,x)^2\, dx+\frac{\beta}{4} \int v(t,x)^4\, dx
$$
Since the energy is conserved we get a global bound for each term in the energy if $\beta\geq 0$ or if $\beta$ is sufficiently small and it follows that in this case the solution constructed in the theorem can be extended to a global solution for $-\infty <t<\infty$. \end{remark}
For the proof we start by introducing the hyperbolic coordinates
$$
\ro^2=t^2-x^2, \quad t=\ro\cosh y,\quad x=\ro\sinh y,
$$
or \beq\label{eq:yweight}
 e^{2|y|}=\frac{t+|x|}{t-|x|},\qquad
\ro^2=t^2-x^2 \eq Then
$$
\square+1=\pa_\ro^2-\ro^{-2}\pa_y^2+\ro^{-1}\pa_\ro+1
$$
and with
$$
v(t,x)=\ro^{-1/2}V(\ro,y)
$$
we get
$$
(\square+1)v(t,x)=\ro^{-1/2}\Big(
\pa_\ro^2+1-\ro^{-2}(\pa_y^2-\tfrac{1}{4})\Big) V(\ro,y).
$$

Hence in these coordinates \eqref{eq:nlkg} becomes the following
equation for $V=\ro^{1/2}v$:
 \beq \label{eq:PDE} \Psi(V)\equiv\pa_\ro^2
V+\Big(1+\frac{\beta}{\ro}V^2+\frac{1}{4\ro^2}\Big)V-\frac{1}{\ro^{2}}\pa_y^2
V=0 \eq
 We are therefore led to first studying the ODE
 \beq\label{eq:ODE}
L(g)\equiv \ddot{g}+\Big(1+\frac{\beta}{\ro}
g^2+\frac{1}{4\ro^2}\Big)g=0 \eq

 We will prove in the next section:
 \begin{prop}\label{prop:firstapproxsolODE} For any constants $a$ and $b$ let
 \beq\label{eq:firstapproxsolODE}
 g_0(\ro)=a\cos\phi,\qquad \phi=\ro+\delta\ln \ro+b,\qquad \delta
 =\frac{3}{8}\beta a^2
 \eq
 Then, if $\delta\geq 0$ is sufficiently small,
 the ODE \eqref{eq:ODE} has a solution $g$ satisfying
 \beq
|\dot{g}-\dot{g}_0|+|g-g_0|\leq C\frac{|a|}{\ro},\qquad \ro\geq 1.
 \eq
 \end{prop}
\begin{remark} The importance of the above proposition is that the ODE
\eqref{eq:ODE} determines the correct phase function $\phi$ in the ansatz for the solution of the PDE \eqref{eq:PDE}.
The precise form of the logarithmic correction to the phase
is due to the long range nature of the interaction.

\end{remark}
\begin{theorem}\label{theorem:asymptotic} Let $v$ be the solution in Theorem
\ref{theorem:exist} and let $V=\rho^{1/2} v$.
Then for each $k\geq 1$ there is $V_k$ of the form
($N$,$I$ finite)
\begin{equation}
a(y) \cos{\phi} + \sum_{n=0}^N\sum_{i\leq  I,\,1\leq j\leq k  }
\Big( a_{i jn}(y) \cos{n\phi} +b_{i jn}(y)
\sin{n\phi}\Big)\frac{\ln^{i}{\!\ro}}{\ro^{j}},
\end{equation}
such that \beq |\Psi(V_k)|\leq \frac{C}{\ro^{k+1}},\qquad |V-V_k|\leq
\frac{C}{\ro^{k+1}}, \qquad \ro \geq 1.
 \eq

 Furthermore, $ a_{ijn}, b_{ijn} $ are monomials in $a$ and its
 derivatives, of at least order $1$ .
\end{theorem}

\section{The first order asymptotics and small data existence at infinity for the ODE}
We want to solve the ODE \eqref{eq:ODE}, subject to a given
behavior at infinity.
\begin{lemma}\label{lemma:approxsolODE}
Let, for any $ |a| \leq 1 $, \beq\label{eq:approxsolODE}
 g_1(\ro)=a\cos\phi+\frac{\delta}{12\ro} a\cos 3\phi,\qquad \phi=\ro+\delta\ln \ro,\qquad \delta
 =\frac{3}{8}\beta a^2
 \eq
 Then
 \beq
|L(g_1)|\leq K\frac{|a|(1+\delta)^2}{\ro^2},\qquad\ro\geq 1
 \eq
 \end{lemma}
\begin{proof}
Using that $\cos^3\phi=\big(\cos
 3\phi+3\cos\phi)/4$ we have $ (A=A(t))$
\begin{equation}
L(A\cos\phi)=
 \Big((1-\dot{\phi}^2)A +\frac{3\beta}{4\ro} A^3+\ddot{A}
 +\frac{A}{4\ro^2}\Big)\cos\phi
 -(\ddot{\phi}A+2 \dot{\phi}\dot{A})\sin\phi
 +\frac{\beta}{4\ro} A^3\cos 3\phi
 \end{equation}
 We get, for $ A(t)=a $ \beq L\big(
a\cos\phi\big)=\Big(-\frac{\delta^2}{\ro^2}+\frac{1}{4\ro^2}\Big)
\,a\cos\phi +\frac{\delta}{\ro^2}\, a\sin\phi
+\frac{2\delta}{3\ro}\, a\cos 3\phi\eq and

 \begin{multline}
L^\prime(0)\Big(-\frac{\cos 3\phi}{8\ro}\Big) =\frac{\cos
3\phi}{\ro} +
\Big(9(\dot{\phi}^2-1)+\frac{1}{\ro}+\frac{1}{4\ro^2}\Big)\frac{\cos
3\phi}{8\ro} +
\Big(\ddot{\phi}-\frac{\dot{\phi}}{\ro}\Big)\frac{3\sin
 3\phi}{8\ro}\\
=\frac{\cos 3\phi}{\ro}+ \frac{\cos
3\phi}{8\ro^2}\big(1+18\delta+\frac{9\delta^2}{\ro}+\frac{1}{4\ro}\big)
-\frac{3\, \sin 3\phi}{8\ro^2}\big(1+\frac{2\delta}{\ro}\big)
 \end{multline}
 where $L^\prime(0)$ is given by \eqref{eq:linODE}.
 Therefore,
 $$ L(g_{1})=L(a\cos\phi) -\frac{8\delta a}{12} L^\prime(0) \big(
 \frac{-\cos 3\phi}{8\ro} \big)
 +\frac{\beta}{\ro}\big(g_1^3-a^3\cos^3\theta\big)= O \big( \frac{(1+\delta)^2 |a|}{\ro^2}
 \big)$$
\end{proof}

 The linearized operator of $L$
around $0$ is given by
 \beq \label{eq:linODE}
L^\prime(0) g=\ddot{g}+g+\frac{1}{4\ro^2}g=F
 \eq
The inverse to this operator with vanishing data at $\infty$ is
given by
 \beq \label{eq:fundsol}
g(\ro)=-\int_\ro^\infty E_s(\ro) F(s)\, ds,
 \eq
 where $E_s(\ro)$ is the forward fundamental solution of
\eqref{eq:linODE}, i.e. $E_s(\ro)$ satisfies $L_0 (E_s)=0$ and
$E_s(s)=0$, $E^{\,\prime}_s(s)=1$. The solution of
\eqref{eq:linODE} satisfies
 \beq
\frac{d}{d\ro}
\big(\dot{g}^2+g^2\big)^{1/2}=\frac{\dot{g}}{\big(\dot{g}^2+g^2\big)^{1/2}}
\Big(F-\frac{g}{4\ro^2}\Big)\geq -
|F|-\frac{\big(\dot{g}^2+g^2\big)^{1/2}}{8\ro^2}.
 \eq
 Multiplying by the integrating factor $e^{-1/(8\ro)}$
 we see that
 \beq\label{eq:linODEest} \big( \dot{g}(\ro)^2+g(\ro)^2\big)^{1/2}\leq
e^{1/(8\ro)}\int_\ro^\infty |F(s)|\, ds\leq 2\int_\ro^\infty
|F(s)|\, ds,\qquad \ro\geq 1
 \eq
 Hence \eqref{eq:fundsol} defines a solution of \eqref{eq:linODE}
 with vanishing data at infinity if the integral above is
 convergent.

 We have
 \beq
L(g)-L(g_1)=L^\prime(0)(g-g_1)+\frac{\beta}{\ro}
G(g_1,g-g_1)(g-g_1), \qquad\text{where}\quad
G(g,h)=\big(3g^2+3gh+h^2\big) \eq  Therefore, to solve (1.8) we
now have to solve the equation
$$
L^\prime(0)(g-g_1)=-\frac{\beta}{\ro} G(g_1,g-g_1)(g-g_1)-L(g_1)
$$
This is done by iteration. We therefore define a sequence $h_k$ :
$$
L^\prime(0)(h_{k+1})=-\frac{\beta}{\ro}
G(g_1,h_k)h_k-L(g_1),\qquad k\geq 0,\qquad h_0=0,
$$
where by Lemma 2.1 (equation 2.2)
 \beq
 |L(g_1)|\leq K\frac{|a|(1+\delta)^2}{\ro^2}
  \eq
 We will inductively assume that
 \beq |h_k|\leq
4K\frac{|a|(1+\delta)^2}{\ro}
 \eq
Then
 \beq
|G(g_1,h_k)|\leq C^\prime |a|^2(1+\delta)^4
 \eq
and by \eqref{eq:linODEest} we have for $\ro\geq 1$,
$$
|h_{k+1}|\leq \int_{\ro}^\infty 2\Big( 4\beta C^\prime
|a|^2(1+\delta)^4 K\frac{|a|(1+\delta)^2}{s^2}
+K\frac{|a|(1+\delta)^2}{s^2}\Big)\, ds\leq
4K\frac{|a|(1+\delta)^2}{\ro}
$$
if $\delta\sim\beta a^2$ is sufficiently small. This shows that we
have a bounded sequence $h_k$, and similarly looking at
differences shows that it converges and hence we get a solution to
the ODE.

\section{The first order asymptotic and small data existence at infinity for the PDE}
In this section we prove Theorem \ref{theorem:exist}
in the case of small data, or equivalently small $\beta$.
This result follows from the general proof in section 4 but
we want to first give the proof in the simple situation were the
complete asymptotic expansion is not needed and one can clearly see
that existence for the PDE follows from existence for the ODE.

We now use Proposition \ref{prop:firstapproxsolODE} and Lemma
\ref{lemma:approxsolODE} to postulate the
following form for the  ansatz of the leading behavior of the
solution of (1.1):
$$
v_1(t,x)=\ro^{-1/2} V_1(\ro,y),
$$
where \beq\label{eq:yapproxsolODE}
 V_1=a(y)\cos\phi(\ro,y)+\frac{\delta}{12\ro} a(y)\cos 3\phi(\ro,y),\qquad
 \phi(\ro,y)=\ro+\delta\ln \ro+b(y),\qquad \delta
 =\frac{3}{8}\beta a^2.
 \eq
Here $a,b$ are smooth functions of $y$, such that $a$ and $b_1=b-b_0$,
where $b_0$ is a constant, are decaying exponentially
 fast. Note that this ansatz is obtained from Lemma 2.1 by simply making the constants $a,b$ dependent on $y$.
  Here we assume that for all $N$,
\beq |D_y^k a(y)|+|D_y^k b_1(y)|\leq C_N e^{-N|y|} \leq C_N
\Big(\frac{t-|x|}{t+|x|}\Big)^{N/2} \leq C_N \frac{\ro^{N}}{t^N}
\eq where we used (1.6) for the second inequality and $ |x| \leq
t$. Therefore, for any $N$,
 \beq
 \frac{|\,a^{(k)}|}{\ro^{N}}\leq \frac{C_N}{t^N}
 \eq

 With notation as in \eqref{eq:ODE} and \eqref{eq:PDE} we have
 \beq
\square \, v_1+v_1+\beta v_1^3=\ro^{-1/2}\Psi(V_1)=\ro^{-1/2}
L(V_1)-\ro^{-1/2} \frac{1}{\ro^2}\pa_y^2 V_1=F_1
 \eq
where if we choose $N$ sufficiently large
 \beq\label{eq:inhomogeneousdecay}
 |F_1|\leq C_N \frac{\big( 1+\beta
 \ln{|1+\ro|}\big)^2}{\ro^{5/2}}\,\,e^{-N|y|}\leq
 C\frac{\big(1+\beta\ln{|1+t|}\big)^2}{t^{5/2}}
 \eq
 since $e^{-2|y|}=(t-|x|)/(t+|x|)$ and $|x| \leq t$.

 We now estimate the correction to $v_1$: let $v$ be the exact
 solution of (1.1).
 We have
\beq (\square +1)(v-v_1)= \beta
G(v_1,v-v_1)(v-v_1)-F_1,\qquad\text{where}\quad
G(v,w)=\big(3v^2+3v w+w^2) \eq
 Let $w$ be the solution of
$$
\square w+w=F
$$
with vanishing data at infinity, i.e. $w$ is defined by
 $$
 w(t,x)=-\int_t^\infty \int E(t-s,x-y) F(s,y)\, dy ds
 $$
 where $E$ is the forward fundamental solution of $\square +1$.
 By the energy inequality
 \beq
 \|\pa w(t,\cdot)\|_{L^2}+\|w(t,\cdot)\|_{L^2}
 \leq \int_t^\infty \|F(s,\cdot)\|_{L^2}\, ds.
 \eq
 Again, we solve for $w=v-v_1$ by iteration:
 Let $w_k$ be defined by $w_0=0$ and
 \beq
(\square+1) w_{k+1}=\beta G(v_1,w_k)w_k+F_1,\qquad k\geq 0.
 \eq
 Since $F_1$ is supported in $|x|\leq t$ it follows from
 \eqref{eq:inhomogeneousdecay} that
 \beq
 \|F_1(t,x)\|_{L^2}\leq \frac{K(1+\beta\ln{|1+t|})^2}{t^{2}}.
 \eq
We will inductively assume that
 \beq \label{eq:inductiveassumption}
 \|\pa w_k(t,\cdot)\|_{L^2}+\|w_k(t,\cdot)\|_{L^2} \leq \frac{4K
(1+\beta\ln{|1+t|})^2}{t}
 \eq
 Since by H\"older's inequality
 $$
 w^2\leq 2\int |w| |w_x|\, dx\leq
 2\|w\|_{L^2} \|\pa w\|_{L^2}\leq \|\pa w\|^2_{L^2}+\|w\|_{L^2}^2
 $$
  we also get
 $$
\|w_k(t,\cdot)\|_{L^\infty}\leq \frac{4K(1+\beta\ln{|1+t|})^2}{t}
 $$
 Since also
(see (3.1) and (3.3))
 \beq
  \| v_1(t,\cdot)\|_{L^\infty}\leq \frac{2C_0}{t^{1/2}}
 \eq
 it follows that for $t \geq t_K $, where $t_K$ depends on $K$ only,
 \beq
 \|G(v_1,w_k)(t,\cdot)\|_{L^\infty} \leq 3\|v_1\|_{L^\infty}^2
 + 3\|v_1\|_{L^\infty} \|w_k\|_{L^\infty} +\|w_k\|_{L^\infty}^2
 \leq \frac{48C_0^2}{t},\qquad t\geq
 t_K
 \eq
Hence by the energy inequality (3.7), and (3.8), (3.9), (3.11)
$$
\|w_{k+1}(t,\cdot)\|_{L^2}\leq \int_t^\infty \big(\beta 8C_0^2+1)
\frac{K(1+\beta\ln{|1+s|})^2}{s^2} ds\leq
\frac{2K(1+\beta\ln{|1+t|})^2}{t},\qquad t\geq t_K^{\,\prime}
$$
if $\beta>0$ is sufficiently small and $t_K^{\,\prime}$ is
sufficiently large. Estimating $\|\pa w_{k+1}\|_{L^2}$ in the same
way, we conclude that \eqref{eq:inductiveassumption} follows also
for $k+1$.

\section{Higher order asymptotics and existence for large data at
infinity}
Let us also consider the linearized operator at $a\cos\phi$:
$$
L_{0}(g)=L^\prime(g_0)\,g=L^\prime(a\cos\phi)
g=\ddot{g}+(1+4^{-1}\ro^{-2}+8\delta\ro^{-1}\cos^2\phi) g
$$
\begin{lemma}
Suppose that $k\geq 1$. We have
\begin{align}
 L_0\Big(\frac{\cos{n\phi}}{\ro^k}\ln^i{\!\!\ro}\Big)&=(1\!-\!n^2)
 \frac{\cos{n\phi}}{\ro^{k}}\ln^i{\!\!\ro}
 +\!\!\!\!\sum_{k^{\,\prime\!}=k+1}^{k+2}\, \sum_{i^{\,\prime}\leq
 i}\,
 \sum_{n^\prime=n-2}^{n+2}\!\!\!\Big(a_{k^{\,\prime\!}
i^{\,\prime}n^{\,\prime}}^{kin}\!
\frac{\cos{n^{\,\prime}\phi}}{\ro^{k^{\,\prime}}}
 +b_{k^{\,\prime\!}
i^{\,\prime}n^{\,\prime}}^{kin}\!\frac{\sin{n^{\,\prime}\phi}}{\ro^{k^{\,\prime}}}\Big)
 \ln^{i^{\,\prime\!}}{\!\!\!\ro}\label{eq:cosn}\\
 L_0\Big(\frac{\sin{n\phi}}{\ro^k}\ln^i{\!\!\ro}\Big)&=(1\!-\!n^2)
 \frac{\sin{n\phi}}{\ro^{k}}\ln^i{\!\!\ro}
 +\!\!\!\!\sum_{k^{\,\prime\!}=k+1}^{k+2}\,\sum_{i^{\,\prime}\leq
 i}\,
 \sum_{n^\prime=n-2}^{n+2}\!\!\!\Big(c_{k^{\,\prime\!}
i^{\,\prime}n^{\,\prime}}^{kin} \!
\frac{\cos{n^{\,\prime}\phi}}{\ro^{k^{\,\prime}}}
 +d_{k^{\,\prime\!}
i^{\,\prime}n^{\,\prime}}^{\,kin}\!\frac{\sin{n^{\,\prime}\phi}}{\ro^{k^{\,\prime}}}\Big)
 \ln^{i^{\,\prime\!}}{\!\!\!\ro}\label{eq:sinn}
 \end{align}
 and
\begin{align}
 L_0\Big(\frac{\cos\phi}{\ro^{k}}\ln^i{\!\!\ro}\Big)&
 =\Big(2k\frac{\sin\phi}{\ro^{k+1}}
 +4\delta \frac{\cos\phi}{\ro^{k+1}}
 +2\delta\frac{\cos 3\phi}{\ro^{k+1}}\Big)\ln^i{\!\!\ro}
 -2i\frac{\sin\phi}{\ro^{k+1}}\ln^{i-1}{\!\!\ro}
 +\!\!\!\sum_{i^{\,\prime\!}\leq i}\!\big(a_{i^{\,\prime\!}}\frac{\cos\phi}{\ro^{k+2}}
 +b_{i^{\,\prime\!}}\frac{\sin\phi}{\ro^{k+2}}\big)\ln^{i^{\,\prime\!}}{\!\!\!\ro}
 \label{eq:cos}\\
 L_0\Big(\frac{\sin\phi}{\ro^{k}}\ln^i{\!\!\ro}\Big)&
 =\Big(-2k\frac{\cos\phi}{\ro^{k+1}}+2\delta\frac{\sin
3\phi}{\ro^{k+1}}\Big)\ln^i{\!\!\ro}
+2i\frac{\cos\phi}{\ro^{k+1}}\ln^{i-1}{\!\!\ro}
+\sum_{i^{\,\prime\!}\leq i}
 \big(c_{i^{\,\prime\!}}\frac{\cos\phi}{\ro^{k+2}}
 +d_{i^{\,\prime\!}}\frac{\sin\phi}{\ro^{k+2}}\big)\ln^{i^{\,\prime\!}}{\!\!\!\ro}
 \label{eq:sin}
 \end{align}
\end{lemma}
\begin{proof}
 Since $\phi=\ro+\delta\ln\ro+b$ it follows that
 \begin{align}
 \frac{d}{d\ro}\Big( \frac{e^{in\phi}}{\ro^k}\Big)&=
\frac{d}{d\ro}\big(e^{in\ro} e^{(i\delta n-k)\ln\ro}e^{ibn}\big)
=\Big(in+\frac{(n\delta
i-k)}{\ro}\Big)\frac{e^{in\phi}}{\ro^k}\\
\frac{d^2}{d\ro^2}\Big( \frac{e^{in\phi}}{\ro^k}\Big)&=
\frac{d^2}{d\ro^2}\big(e^{in\ro} e^{(i\delta
n-k)\ln\ro}e^{ibn}\big) =\Big(-n^2+\frac{2ni(n\delta
i-k)}{\ro}+\frac{c_{nk}}{\ro^2}\Big)\frac{e^{in\phi}}{\ro^k}
 \end{align}
Hence \beq \frac{d^2}{d\ro^2}\Big(
\ln^i{\!\!\ro}\frac{e^{in\phi}}{\ro^k}\Big)
=\Big(-n^2+\frac{2ni(n\delta
i-k)}{\ro}+\frac{c_{nk}}{\ro^2}\Big)\frac{e^{in\phi}}{\ro^k}
\ln^i{\!\!
\ro}+\sum_{k^{\,\prime}=k+1}^{k+2}\sum_{i^{\,\prime}=i-2}^{i-1}c_{kik^{\,\prime}
i^{\,\prime}} \frac{e^{in\phi}}{\ro^{k^{\,\prime}}}
\ln^{i^{\,\prime}}{\!\!\ro} \eq

\begin{align}
\frac{d^2}{d\rho^2}\Big(\frac{\cos{n\phi}}{\ro^k}\Big)
&=-n^2\frac{\cos{n\phi}}{\ro^k}-2\delta
n^2\frac{\cos{n\phi}}{\ro^{k+1}} +2 k n
\frac{\sin{n\phi}}{\ro^{k+1}} +a_{kn}\frac{\cos{n\phi}}{\ro^{k+2}}
 +b_{kn}\frac{\sin{n\phi}}{\ro^{k+2}}\label{eq:twodercos}\\
\frac{d^2}{d\rho^2}\Big(\frac{\sin{n\phi}}{\ro^k}\Big)
&=-n^2\frac{\sin{n\phi}}{\ro^k}-2\delta
n^2\frac{\sin{n\phi}}{\ro^{k+1}} -2 k n
\frac{\cos{n\phi}}{\ro^{k+1}} +c_{kn}\frac{\cos{n\phi}}{\ro^{k+2}}
 +d_{kn}\frac{\sin{n\phi}}{\ro^{k+2}}\label{eq:twodersin}
\end{align}

 Since
$\cos^2\phi\cos{\phi}=\big(3\cos\phi+\cos 3\phi\big)/4$ and
$\cos^2\phi\sin{\phi}=\big(\sin\phi+\sin 3\phi\big)/4$, we have
\begin{align}
\Big(1+8\delta\frac{\cos^2\phi}{\ro}\Big)\frac{\cos{\phi}}{\ro^k}
&=\frac{\cos\phi}{\ro^k}+6\delta\frac{\cos\phi}{\ro^{k+1}}+2\delta\frac{\cos
3\phi}{\ro^{k+1}}\\
\Big(1+8\delta\frac{\cos^2\phi}{\ro}\Big)\frac{\sin{\phi}}{\ro^k}
&=\frac{\sin\phi}{\ro^{k}}+2\delta\frac{\sin\phi}{\ro^{k+1}}+2\delta\frac{\sin
3\phi}{\ro^{k+1}}
\end{align}
\end{proof}

\begin{defi} Let ${\cal S}_k$ denote the family of finite sums ($N$,$I$,$j$-sum finite) of the form
\begin{equation}
\sum_{n=0}^N\sum_{i\leq  I,\, j\geq k  } \Big( a_{i jn}(y)
\cos{n\phi} +b_{i jn}(y)
\sin{n\phi}\Big)\frac{\ln^{i}{\!\ro}}{\ro^{j}},\qquad
\phi=\ro+\delta\ln\ro,\quad \delta=\frac{3}{8}\beta a(y)^2
\end{equation}
where for any $N$ and $\ell$ there is a constant such that
\begin{equation}
|D_y^\ell a_{ijk}(y)|+|D_y^\ell b_{ijk}(y)|\leq C_{N\ell\, } e^{-N|y|}
\end{equation}
Furthermore, let $\ringabove{\cal S}_k$ denote the family of finite sums of
the above form but with
 $$
 a_{ik1}=b_{ik1}=0,\qquad\text{for all}\quad i
 $$
\end{defi}
\begin{lemma}\label{lemma:L0inv}
If $k\geq 1$ and $\ringabove{\Sigma}_k\in \ringabove{\cal
S}_k$, then there are $\Sigma_k\in {\cal S}_k$ and
$\ringabove{\Sigma}_{k+1}\in\ringabove{\cal S}_{k+1}$ such that
\begin{equation}
L_0 \Sigma_k =\ringabove{\Sigma}_k +\ringabove{\Sigma}_{k+1}
\end{equation}
\end{lemma}
\begin{proof} First we use the first part of the previous lemma to
invert the terms with $k^{\,\prime\!}=k$ and $n\neq 1$. Then we
use the second part of the previous lemma to successively remove
the terms with $n=1$ by lowering the logarithms. First note that
an element of $ \ringabove {\cal S}_k$ can be written as a
$$ \sum_{n\neq 1} \sum_{k^\prime \geq k,\, i} \big(\alpha_{ik^\prime  n}
\frac {\cos n\phi}{\ro^{k^\prime}} +\beta_{ik^\prime n}
\frac{\sin n\phi}{\ro^{k^\prime}}\big)\ln^i(\ro) + \sum_{k^\prime
\geq k+1,\, i}\big( \alpha^\prime_{ik^\prime }
\frac{\cos\phi}{\ro^{k^\prime}} + \beta_{ik^\prime }^\prime
\frac{\sin\phi}{\ro^{k^\prime}}\big) \ln^i(\ro)= I_{nr}+I_{res}$$
The sum over $ k^\prime \geq k+1 $ is due to the fact that $
\ringabove{\Sigma_k} $ unlike $ \Sigma_k $ is "nonresonant", that
is do not contain lowest order terms in $\ro$ for $ n=1$. (see
definition of the space $ \ringabove{S_k} $). Now, given such
element $ \ringabove{\Sigma_k} $, we use (4.1), (4.2) to obtain
$$ L_0 \big( \sum_{n\neq 1} \frac{1}{1-n^2} \sum_{k^\prime \geq
k,\, i} \big ( \alpha_{ik^\prime n}\cos n\phi + \beta_{ik^\prime
n}\sin
 n\phi \big) \frac {\ln^i\ro}{\ro^{k^\prime}} \big)= I_{nr} +
\ringabove{\Sigma_{k+1} }
$$
We are therefore left with inverting $L_0$ on $I_{res}$. To this
end we use (4.3),(4.4), to obtain :
\begin{multline}
 L_0 \big( \frac{1}{2k} \frac{\cos \phi}{\ro^k} \ln^i\ro
+ \frac{4\delta}{2k} \frac{\sin \phi}{\ro^k} \ln^i\ro \big)=
\frac{\sin \phi}{\ro^{k+1}} \ln^i\ro + \\
O(\ro^{-k-1}\ln^i\ro) (\sin 3\phi , \cos 3\phi) +
O(\ro^{-k-1}\ln^{i-1}\ro) (\sin\phi , \cos\phi) +
O(\ro^{-k-2}\ln^{i}\ro)(\sin\phi , \cos\phi)
\end{multline}
 and similar formula for $ \cos\phi\, \ro^{-k-1}
\ln^i\ro$. (Here $(f,g) \equiv \alpha a + \beta g $ for some
numbers $ \alpha, \beta$.) Hence, we can invert $L_0$ on
$(\cos\phi, \sin\phi)\, \ro^{-k-1} \ln^i\ro$  up to nonresonant
terms in $\ringabove{S_{k+1}} $, $ S_{k+2}$ and resonant terms in
$S_{k+1}$ but with one less power of $\ln\ro$. Hence, by
iteration, eliminate all such terms, in each step one less power
of $\ln\ro$.
\end{proof}

We must then show that the products of the above classes are
properly mapped as well as the Laplacian acting on the above
classes. We want to solve
\begin{equation}
\Psi(V)=\pa_\ro^2
V+\Big(1+\frac{\beta}{\ro}V^2+\frac{1}{4\ro^2}\Big)V-\ro^{-2}\pa_y^2
V=0
\end{equation}
by iteration, starting from
\begin{equation}
V_0 = a\cos\phi,\qquad\qquad \text{where}\quad \phi=\ro+3
\beta 8^{-1} a^2\ln \ro+b
\end{equation}
 and $a=a(y)$, $b=b(y)$.
 We have:
 \begin{lemma} There is a sequence $V_k$, $k=0,...$, such $V_k-V_0\in {\cal S}_1$,
$\Psi(V_k)\in \ringabove{\cal S}_{k+1}$ and $V_{k}-V_{k-1}\in {\cal S}_k$.
\end{lemma}
 \begin{proof}We have
 \beq
 \pa_y^2 \big(a\cos\phi)=\sum_{j=0}^2 \big(a_j
 \cos\phi+b_j\sin\phi\big)\ln^{j\!\!}\ro
 \eq
 for some functions $a_j(y)$, and $b_j(y)$ which are at least
 linear in $ a^{(k)}, b^{(k)}$, $ (0 \leq k \leq 2)$.
 It follows that
\begin{equation}
\Psi(V_0)=-\frac{\delta^2}{\ro^2} \,a\cos\phi
+\frac{\delta}{\ro^2}\, a\sin\phi +\frac{2\delta}{3\ro}\, a\cos
3\phi+\frac{a\cos\phi}{4\ro^2}-\frac{1}{\ro^2}\pa_y^2
\big(a\cos\phi)\in \ringabove{\cal S}_1
\end{equation}
This proves the lemma for $k=0$ and in what follows we will assume the lemma
for $k$ replaced by $k-1$ and show that this implies the lemma also for $k$.

We have
\begin{equation}
\Psi^\prime(V)W=\pa_\ro^2
W+\Big(1+3\frac{\beta}{\ro}V^2+\frac{1}{4\ro^2}\Big)W-\ro^{-2}\pa_y^2
W
\end{equation}
 Since the operator
$\ro^{-2}\pa_y^2\psi$ maps ${\cal S}_k\to {\cal
S}_{k+2}\subset \ringabove{\cal S}_{k+1}$ it follows that
\begin{equation}
\Psi^\prime(V_0)= L_0 - \ro^{-2} \pa_y^2 V
\end{equation}
can be inverted in the same spaces as $L_0$ in Lemma
\ref{lemma:L0inv}, i.e.
if $k\geq 1$ and $\ringabove{\Sigma}_k\in \ringabove{\cal
S}_k$, then there are $\Sigma_k\in {\cal S}_k$ and
$\ringabove{\Sigma}_{k+1}\in\ringabove{\cal S}_{k+1}$ such that
 \beq
\Psi^{\,\prime}(a\cos\phi)\Sigma_k=\ringabove{\Sigma}_{k}+\ringabove{\Sigma}_{k+1}
 \eq
 Moreover if $V_n-V_0\in {\cal S}_1$ it follows that
 $\big(\Psi^\prime(V_n)-\Psi^\prime(V_0)\Big)\Sigma_k=3\beta (V_n^2-V_0^2)\Sigma_k/\ro\in {\cal S}_{k+2}\in \ringabove{\cal S}_{k+1}$ so
 $\Psi^\prime(V_n)$ also satisfies
 \beq\label{eq:psikinv}
\Psi^{\,\prime}(V_n)\Sigma_k=\ringabove{\Sigma}_{k}+\ringabove{\Sigma}_{k+1}
 \eq
 for some other $\ringabove{\Sigma}_{k+1}$.

Given $V_{k-1}$ such that $\Psi(V_{k-1})\in \ringabove{\cal S}_{k}$
and $V_{k-1}-V_0\in {\cal S}_1$ we now
find $V_{k}$ such that $V_{k}-V_{k-1}\in {\cal S}
_{k}$ by solving
 \beq
\Psi^{\,\prime}(V_{k-1})(V_{k}-V_{k-1})+\Psi(V_{k-1})\in \ringabove{\cal
S}_{k+1}
 \eq
 which is possible, by \eqref{eq:psikinv}.
 Then with $\Phi(V,U)=3 V+U$
\beq \Psi(V_{k})=\Psi(V_{k-1})+\Psi^{\,\prime}(V_{k-1})(V_{k}-V_{k-1})
+\frac{\beta}{\ro}\Phi(V_{k-1},V_{k}-V_{k-1})(V_{k}-V_{k-1})^2\in
\ringabove{\cal S}_{k+1}
 \eq
\end{proof}

We have now found $v_N$, for any $N$ , such that
$$
\square v_N+v_N+\beta v_N^3=F_N=O(t^{-N-5/2}),\qquad
v_N-v_0=O(t^{-3/2}\ln{t})
$$
It follows that there is a constant $C_0<\infty$ independent of
$N$ and another constant $t_N<\infty$ depending on $N$ such that
$$
|v_N|\leq 2C_0 t^{-1/2},\qquad t\ge t_N
$$
We then define $w_0=0$ and for $l\geq 1$:
$$
(\square+1) w_{l+1}=\beta G(v_{N},w_l)w_l+F_N,\qquad l\geq 0.
 $$
 Since $F_N$ is supported in $|x|\leq t$ it follows from
 \eqref{eq:inhomogeneousdecay} that
 $$
 \|F_N(t,x)\|_{L^2}\leq \frac{K_N}{t^{N+1}}
 $$
We will inductively (in $l$)  assume that
 \beq\label{eq:inductiveassumptionN}
 \|\pa w_l(t,\cdot)\|_{L^2}+\|w_l(t,\cdot)\|_{L^2} \leq \frac{4K_N}{Nt^{N}}
 \eq
 Since by H\"older's inequality
 $$
 w^2\leq 2\int |w| |w_x|\, dx\leq
 2\|w\|_{L^2} \|\pa w\|_{L^2}\leq \|\pa w\|^2_{L^2}+\|w\|_{L^2}^2
 $$
  we also get
 $$
\|w_l(t,\cdot)\|_{L^\infty}\leq \frac{4K_N}{N\,t^N}
 $$
 Since also
 $$
\| v_N(t,\cdot)\|_{L^\infty}\leq \frac{2C_0}{t^{1/2}},\qquad t\geq
t_N
 $$
 where $C_0$ is independent of $N$,
 it follows that
 $$
 \|G(v_N,w_l)(t,\cdot)\|_{L^\infty}\leq \frac{8C_0}{t},\qquad t\geq
 t_N^{\,\prime}
 $$
Hence by the energy inequality
$$
\|\pa w_{l+1}(t,\cdot)\|_{L^2}+\|w_{l+1}(t,\cdot)\|_{L^2}\leq
\int_t^\infty \frac{\beta 8C_0}{s} \frac{ 4K_N}{N s^N}
ds+\frac{K_N}{s^{N+1}}\, ds=\Big(\frac{32\,\beta
C_0}{N}+1\Big)\frac{K_N}{N\, t^N} \leq \frac{2K_N}{N t^N} ,\qquad
t\geq t_N^{\,\prime\prime}
$$
if  $t_N^{\,\prime\prime}$ and $N$ are sufficiently large. Hence
\eqref{eq:inductiveassumptionN} follows also for $l+1$.


\begin{thebibliography}{C-K}
\expandafter\ifx\csname url\endcsname\relax
  \def\url#1{{\tt #1}}\fi
\expandafter\ifx\csname
urlprefix\endcsname\relax\def\urlprefix{URL }\fi



\bibitem[D1]{D1}
J.-M. Delort
\newblock {\em Existence globale et comportement asymptotique pour
l'\'equation de Klein-Gordon quasi lin\'eaire \`a donn\'ees petites en
dimension 1. \/.}
\newblock Ann. Sci. \'Ecole Norm. Sup.  no. 4,
\textbf{34} (2001), 1--61.

\bibitem[H1]{H1} L. H\"ormander \newblock
Lectures on Nonlinear hyperbolic differential equations.
\newblock Springer Verlag (1997)


\bibitem[L1]{L1} H. Lindblad,
\newblock {\em Global solutions of nonlinear wave equations\/.}
\newblock Comm.Pure Appl. Math. \textbf{45} (9) (1992), 1063-1096.


\bibitem[L-R]{L-R}
H. Lindblad and I. Rodnianski
\newblock {\em The weak null condition for Einstein's equations\/.}
\newblock C. R. Math. Acad. Sci. Paris 336 (2003), no. 11,
901--906

\end{thebibliography}
\end{document}